**Х. К. Дакхіл, Ю. Б. Зелінський, Б. А. Кліщук**

КНУ ім. Т.Г. Шевченка, Інститут математики НАН України, Київ

*E-mail*: moon5385@gmail.com, zel@imath.kiev.ua, bogdanklishchuk@mail.ru

**Про слабко *m*–опуклі множини**

*Досліджено деякі властивості слабко m-опуклих множин. Отримано оцінки для різних варіантів задачі про тінь в фіксованій точці. Побудовано приклад, що дає оцінку знизу необхідної кількості куль для створення тіні в точці, дотичної до сфери $S^2$ в просторі $\mathbb{R}^3$*

**Ключові слова:** *m-опукла множина, слабко m-опукла множина, грасманів многовид, спряжена множина, задача про тінь, 1-оболонка сім'ї множин*

В даній статті розглянуто різни варіанти задачі про тінь в точці в просторі $\mathbb{R}^n$, зокрема досліджено задачу про тінь в точці, дотичну до сфери $S^2 \subset \mathbb{R}^3$.

Далі під *m*-вимірними площинами розуміємо *m*-вимірні афінні підпростори евклідового простору $\mathbb{R}^n$.

**Означення 1.** Скажемо, що множина $E \subset \mathbb{R}^n$ *m*-опукла відносно точки $x \in \mathbb{R}^n \setminus E$, якщо знайдеться *m*-вимірна площина L, така що $x \in L$ і $L \cap E = \varnothing$.

**Означення 2.** Скажемо, що відкрита множина $G \subset \mathbb{R}^n$ слабко *m*-опукла, якщо вона *m*-опукла відносно кожної точки $x \in \partial G$, яка належить межі множини G. Скажемо, що довільна множина $E \subset \mathbb{R}^n$ слабко *m*-опукла, якщо її можна аппроксимувати ззовні сім'єю відкритих слабко *m*-опуклих множин.

Легко побудувати приклад слабко *m*-опуклої, але не *m*-опуклої множини.

**Приклад 1.** Нехай $B = \{(x,y) | x^2 + y^2 < 1\}$ – відкритий круг на площині *xOy*. Виберемо три точки на колі $S^1 = \{(x,y) | x^2 + y^2 = 1\}$ і розглянемо симплекс σ з вершинами в цих точках. Неважко переконатися, що множина $E = B \setminus \sigma$ слабко 1-опукла, але не 1-опукла.

В роботах [1-3] досліджено деякі властивості слабко *m*-опуклих та *m*-опуклих множин.

Нехай $G(n,m)$ – грасманів многовид *m*–вимірних площин в $\mathbb{R}^n$ [4]. Спряженою множиною $E^*$ до множини $E$ називається підмножина *m*-вимірних площин в $G(n,m)$, які не перетинають множину $E$.

Дослідимо деякі властивості слабко *m*-опуклих множин.

**Твердження 1.** *Якщо $E_1$ і $E_2$ відповідно слабко k-опукла і слабко m-опукла множина $k \leq m$, то $E = E_1 \cap E_2$ буде слабко k-опуклою множиною.*

**Теорема 1.** *Якщо K – слабко m-опуклий компакт і множина $K^*$ зв'язна, то для перерізу K довільною (n − m)-вимірною площиною L множина $L \setminus K \cap L$ зв'язна.*

*Доведення.* Як встановлено в твердженні 2 [2] множина $K^*$ є відкритою множиною, тому будь які дві її точки можна з'єднати неперервною дугою в $K^*$. Припустимо, що існує така (n − m)-вимірна площина *L*, для якої множина $L \setminus K \cap L$ незв'язна. Отже переріз є носієм деякого ненульового (n − m − 1)-мірного циклу *z* [4]. Нехай точка *x* належить обмеженій компоненті множини $L \setminus K \cap L$. Такі точки існують внаслідок компактності *K*. Не порушуючи загальності, внаслідок слабкої *m*-опуклості *K*, можна вважати, що через точку *x* можна провести *m*-вимірну площину $l_1$, яка не перетинає *K*. Іншу *m*-вимірну площину $l_2$ проведемо за межами деякої досить великої кулі, що містить компакт *K*. Якщо ми компактифікуємо простір $\mathbb{R}^n$ до сфери $S^n$ нескінченно віддаленою точкою, то отримаємо два *m*-мірні цикли $w_1 = l_1 \cup (\infty)$ і $w_2 = l_2 \cup (\infty)$, з яких перший цикл зачеплений з циклом *z*, а другий ні. З одного боку, ці цикли не можна перевести гомотопією один в інший, яка не перетинала б цикл *z*, а, отже і множину $K \cap L$. З другого боку, зв'язність множини $K^*$ забезпечує існування в $K^*$ пари точок $y_1$, $y_2$, які задають площини $l_1$ і $l_2$ відповідно і з'єднані дугою в $K^*$. Точки цієї дуги задають гомотопію площини $l_1$ в $l_2$, яка не має спільних точок з множиною $K \cap L$. Отримане протиріччя доводить теорему.

В роботі [5] досліджено різні варіанти наступної задачі, яку популярно можна назвати задачею про тінь в точці.

**Задача**. Яка мінімальна кількість опуклих множин, які попарно не перетинаються і з деякими заданими умовами, достатня для того, щоби будь яка пряма, що проходить через наперед задану точку, перетинала хоча б одну з цих множин?

Якщо цього досягти, то кажуть, що вибрана точка належить до 1-оболонки даної сім'ї множин.

Для набору куль різного радіуса в [5] отримана оцінка.

**Теорема 2**. *Для того щоб вибрана точка в n-вимірному евклідовому просторі при n ≥ 2 належала 1-оболонці сім'ї відкритих (замкнутих) куль, які дану точку не містять і попарно не перетинаються, необхідно і достатньо n куль.*

Виявляється, що у випадку куль однакового радіуса результат відрізняється від випадку, коли кулі різних радіусів.

**Теорема 3**. *Довільний набір із трьох куль однакового радіуса, які попарно не перетинаються, утворює слабко 1-опуклу множину в тривимірному евклідовому просторі.*

Для набору з трьох куль в просторі $\mathbb{R}^n$ має місце твердження.

**Теорема 4**. *Для довільної точки простору $\mathbb{R}^n \setminus \bigcup_{i=1}^{3} B_i$, де $B_1$, $B_2$, $B_3$ — набір з трьох куль однакового радіуса, які попарно не перетинаються і не проходять через цю точку, існує (n − 2)-вимірна площина, що містить цю точку і не перетинає жодну з куль.*

**Означення 3.** Скажемо, що сім'я множин $\Im = \{F_\alpha\}$ задає тінь, дотичну до многовиду $M$ в точці $x \in M$, якщо кожна пряма, дотична до многовиду $M$ в точці $x \in M \setminus \cup_\alpha F_\alpha$, має непорожній перетин хоча б з однією з множин $F_\alpha$, яка належить до сім'ї $\Im$.

**Лема.** *Розглянемо рівносторонній трикутник ABC в евклідовій площині $R^2$. Якщо ми виберемо три круги $B_i, i = 1,2,3$, з центрами в вершинах цього трикутника і радіуса рівного половині висоти трикутника, то кожна пряма, що проходить через довільну точку $x \in (\bigcup_{i=1}^{3} B_i)^* \setminus \bigcup_{i=1}^{3} B_i$, де $(\bigcup_{i=1}^{3} B_i)^*$ – опукла оболонка множини $\bigcup_{i=1}^{3} B_i$, перетинається не менше ніж з одним із вибраних кругів.*

Відмітимо, що коло, описане навколо цього трикутника, лежить в опуклій оболонці цих трьох кругів. Цей результат показує, що в тривимірному випадку для довільної точки сфери можна вибрати три кулі, які попарно дотикаються і які забезпечать тінь в усіх точках криволінійного трикутника, вирізаного на сфері цими кулями. Але, як показує наступний приклад, узгодження такої конструкції для всієї сфери вимагає додаткових міркувань.

**Приклад 2.** Існує набір із 14 відкритих (замкнутих) куль, що попарно не перетинаються с центрами на сфері $S^2 \subset \mathbb{R}^3$, який не може забезпечити тінь, дотичну до сфери $S^2$ в кожній точці $x \in S^2 \setminus \bigcup_{i=1}^{14} B_i$.

Не порушуючи загальності, можемо припустити, що вибрана сфера $S^2$ з центром в початку координат радіуса одиниця. Впишемо в цю сферу куб з вершинами в точках $(x = \pm 1/\sqrt{3}, y = \pm 1/\sqrt{3}, z = \pm 1/\sqrt{3})$. Довжина ребра куба дорівнює $a = 2/\sqrt{3}$. Тепер виберемо вісім відкритих куль з центрами в вершинах куба і радіуса рівного половині ребра куба $r = 1/\sqrt{3} \approx 0.577$. Добавимо до цього набору шість нових відкритих куль з центрами в точках перетину променів, які виходять з початку координат і проходять через центр грані куба, зі сферою $S^2$. Радіуси цих куль дорівнюють $r_1 = \sqrt{2 - 2/\sqrt{3}} - 1/\sqrt{3}$. Кожен з них доторкається рівно до чотирьох раніше вибраних куль. Цей набір куль двох різних радіусів покриває сферу, але як показують обчислення, цієї системи куль не вистачає для створення тіні, дотичної до сфери $S^2$ в кожній точці $x \in S^2 \setminus \bigcup_{i=1}^{14} B_i$.

Побудований набір куль дає оцінку знизу необхідної кількості куль. Оцінка зверху залишається відкритою.

*Х. К. Дакхил, Ю. Б. Зелинский, Б. А. Клищук*
КНУ им. Т.Г. Шевченка, Институт математики НАН Украины, Киев
*E-mail*: moon5385@gmail.com, zel@imath.kiev.ua, bogdanklishchuk@mail.ru

**О слабо *m*-выпуклых множествах**

*Исследованы некоторые свойства слабо m-выпуклых множеств. Получены оценки для различных вариантов задачи о тени в фиксированной точке. Построен пример, который дает оценку снизу необходимого количества шаров для создания тени в точке, касательной к сфере $S^2$ в пространстве $\mathbb{R}^3$*

**Ключевые слова:** *m-выпуклое множество, слабо m- выпуклое множество, грассманово многообразие, сопряженное множество, задача о тени, 1-оболочка семейства множеств*

*H. K. Dakhil, Yu. B. Zelinskii, B. A. Klishchuk*
Taras Shevchenko National University of Kyiv, Institute of Mathematics NAS of Ukraine, Kyiv
*E-mail*: moon5385@gmail.com, zel@imath.kiev.ua, bogdanklishchuk@mail.ru

**On weakly *m*-convex sets**

*There were obtained some properties of weakly m-convex sets. Various variants of the problem of shadow were investigated. There was obtained the lower estimation for the number of balls that are necessary to create a shadow at the point of the sphere $S^2$ in Euclidean space $\mathbb{R}^3$*

**Key words:** *m-convex set, weakly m-convex set, Grassmann manifold, conjugate set, problem of shadow, 1-hull of a family of sets*

Дані про авторів

**Хайджаа Кхудхаір Дакхіл**
КНУ ім. Т.Г. Шевченка
*E-mail*: moon5385@gmail.com

**Юрій Борисович Зелінський**
Інститут математики НАН України
*E-mail*: zel@imath.kiev.ua

**Богдан Анатолійович Кліщук**
Інститут математики НАН України
*E-mail*: bogdanklishchuk@mail.ru